\documentclass[12pt]{article}
\usepackage{hyperref,url,graphicx}
\usepackage{xcolor}
\usepackage[version=3]{mhchem}
\usepackage{graphicx}
\usepackage{amssymb,amsmath}
\usepackage[T1]{fontenc}
\usepackage[utf8]{inputenc}

\newtheorem{Thm}{Theorem}
\newtheorem{Def}{Definition}
\newtheorem{Exa}{Example}

\newtheorem{Sta}{Statement}
\newcommand{\Proof}[1]{{\textit{Proof. }#1}}

\newcommand{\torol}[1]{}

\newcommand{\mth}{$m^{\mbox{th}}$\ }
\newcommand{\rth}{$r^{\mbox{th}}$\ }

\newcommand{\alphab}{\boldsymbol{\alpha}}
\newcommand{\betab}{\boldsymbol{\beta}}
\newcommand{\gammab}{\boldsymbol{\gamma}}

\newcommand{\nulb}{\boldsymbol{0}}
\newcommand{\oneb}{\boldsymbol{1}}

\newcommand{\ab}{\mathbf{a}}

\newcommand{\cb}{\mathbf{c}}

\newcommand{\kb}{\mathbf{k}}
\newcommand{\Kb}{\mathbf{K}}
\newcommand{\vb}{\mathbf{v}}
\newcommand{\yb}{\mathbf{y}}

\newcommand{\N}{\mathbb{N}}

\newcommand{\emp}[1]{\textbf{#1}}
\newcommand{\mma}[1]{\textbf{\tt #1}}

\title{Automatic model generation\footnote{Dedicated to Professors L\'aszl\'o Lov\'asz and J\'ozsef P\'alink\'as.}}

\author{
Tibor Nagy$^{1,2}$\ ,
J\'anos T\'oth$^{2,3}$\ ,
Tam\'as Ladics$^{4}$\\
${}^{1}$Institute of Materials and Environmental Chemistry,\\
Research Centre for Natural Sciences,\\
Hungarian Academy of Sciences,
Budapest, Hungary\\
${}^{2}$Laboratory for Chemical Kinetics,\\
E\"otv\"os Lor\'and University,
Budapest, Hungary\\
${}^{3}$Department of Mathematical Analysis\\
Budapest University of Technology and Economics,\\
Budapest, Hungary\\
${}^{4}$Department of Science and Engineering,\\
John von Neumann University,
Kecskem\'et, Hungary\\
Correspondence
Tibor Nagy,\\
Institute of Material and Environmental Chemistry,\\
Research Centre for Natural Sciences,\\
Hungarian Academy of Sciences,\\
1117 Budapest, Magyar tud\'osok k\"or\'utja 2, Hungary\\
Email : \url{nagy.tibor@ttk.mta.hu}\\
Funding information\\
National Research Development and Innovation Office,\\ SNN 125739.}

\begin{document}
\maketitle
\begin{abstract}
The goal of the paper is to automatize the selection of mechanisms which are able to describe
a set of measurements.
In order to do so first we construct a set of possible mechanism fulfilling
chemically reasonable requirements with a given number of species and reaction steps.
Then we try to fit all the  mechanisms, and offer the best fitting one to the chemist for further analysis.
The method can also be used to a kind of lumping: to reproduce the results of a big mechanism
using a smaller one, with less number of species.
We show two applications: one on an artificial example and another one on a small real life data.
\end{abstract}
KEYWORDS reaction networks, inverse problem of reaction kinetics,  automatically generated models,  reduced models, model parametrization, fitting rate coefficients
\section{INTRODUCTION}
Our goal is to help solve an inverse problem of reaction kinetics:
we try to build "reaction mechanism templates", candidate set of reactions which are able to serve as models for given (measured or simulated with a detailed model) concentration vs. time curves.
This we do by constructing a large set of reaction mechanisms then by discarding those which do not fulfil some chemically relevant restrictions.
To put it in another way, we provide a set of possible models automatically for the chemist, and say,
which of them fits best to the measurements or the simulation measurements of a detailed model.
The approach of the Ref.
\cite{lokbrent} is quite similar, but (almost) all the details are different.
Their goal is to construct possible cellular reaction networks.

Thus, given are: concentrations as functions of time, and our task is to
find a \emph{model}, a \emph{reaction mechanism}, i.e. a set of reaction steps (also called
complex chemical reaction, reaction network, or simply reaction), endowed with mass action type kinetics and with appropriate reaction rate coefficients from the set of \emph{all} possible models.
Concentrations used as reference are either determined from measurements or from simulations of a detailed model.
Here we show how the method works by testing mainly on \emph{simulated} data, and also a set of real experimental data: measurements on the salicylic acid transport under different conditions..

A related problem is to find \emph{good initial estimates} of the parameters$^{\cite{ladicstothnagy}}$,
in order to accelerate parameter fitting, which needs to be  done for a large number of candidate models.

The structure of our paper is as follows. Section 2 offers a possible scenario to construct candidate sets of reactions which are to describe the concentration vs. time functions.
Section 3 shows the applicability of the method on two examples.
In the first one, a specific model has been chosen from the set of models constructed according to Section 2;
reference data are generated by its simulation.
In the second one we use experimental data on salicylic acid transport from the stomach to the intestine with or without addenda measured at different pH values.
Beyond the original model (the classical irreversible consecutive reaction)
other mechanisms are also tried.
The last Section is about possible extensions and formulates some open problems, mainly in combinatorics.
\section{A SET OF CHEMICALLY REASONABLE RESTRICTIONS}
Basic concepts of formal reaction kinetics are presented here shortly, for the formal and detailed expansion we propose the use of e.g. Refs. \cite{tothnagypapp} or \cite{feinbergbook}.
\subsection{Species}
First, we have to fix the number of species \(M\).
In our illustrating examples it will usually be 2 or 3; in applications this may be equal to the number of measured concentrations.
\subsection{Complexes}
In formal reaction kinetics the linear combinations
of species with stoichiometric coefficients
on the sides of the reaction arrows are called \emph{complexes}, a slightly unfortunate name, because this word is used in chemistry with a completely different meaning. Their number is usually denoted by $N.$ The stoichiometric coefficients of the \mth species on the left side of the \rth reaction arrow is denoted by $\alpha_{mr},$ those on the right side by $\beta_{mr},$ thus the general form of our mechanism is
\begin{equation}\label{eq:R1}
\sum_{m=1}^M\alpha_{mr}\ce{X($m$)}\rightleftharpoons\sum_{m=1}^M\beta_{mr}\ce{X($m$)}.
\end{equation}
\paragraph{Mass conservation}
In systems where chemical reactions significantly modify the mole fraction of species
(e.g. the changes are larger than \(1\%\)), mass
conservation has to be strictly fulfilled.
This is the case in combustion models, and we shall accept this view here.
Whereas atmospheric chemistry models usually contain steps from which major components of air (e.g. \ce{N2}, \ce{O2}, \ce{CO2} in tropospheric chemistry models) are simply omitted if formed as negligible change in their concentration is induced by the step.
The number of steps is less if only mass conserving reactions are allowed.

As a consequence of mass conservation $M=1$ is immediately excluded, because e.g. the reaction step
\ce{X <=> 2X} is not allowed.
Also excluded is the empty complex 0, as it immediately leads to mass destruction or mass creation.
Still, the empty complex may either be useful to describe the sticking of a species
to the wall or leave the reactor in any other way,
or to express an inflow or outflow in formal reactions, e.g. \ce{0 -> X} or \ce{X + Y -> 0}.
\paragraph{Short complexes}
In most cases reaction steps containing complexes longer than 2, are not accepted; here we accept this view.
In other words, we only allow \emph{short complexes}, complexes  of the form in which no more than two species take part. In the case when \(M=2\) we have thus
\ce{X},\ce{Y},\ce{2X},\ce{2Y},\ce{X + Y} (zero complex excluded from the beginning).
Their number in the general case, as it can be immediately seen, is
\begin{equation}\label{eq:nocomplexes}
  \overline{N}(M)=M+M+\binom{M}{2}=\frac{M(M+3)}{2}.
\end{equation}
Let us mention in passing that the concept of short complexes proved to be really useful when systematically investigating the dynamic behaviour of small reaction mechanisms$^{\cite{hornprs3}}$.
However, our present restriction does not mean that termolecular reactions are unimportant.
In gas phase reactions like
\begin{equation}\label{eq:R2}
\ce{2 NO + O2 -> 2 NO2}
\end{equation}
are quite common, see e.g.
\url{https://engineering.columbia.edu/news/michael-burke-termolecular-reactions}.
To take into consideration these one should also include \ce{3 X},\ce{2 X + Y},\dots,\ce{X + Y + Z},
the total number of which is
\[
M+M(M-1)+\binom{M}{3}=\frac{M(M+1)(M+2)}{6}=\binom{M+2}{3}.
\]
\subsection{Reaction steps}
\paragraph{Reversibility}
All elementary reactions are strictly reversible due to microscopic reversibility.
However, in some cases the rate of the forward or backward reaction
can be negligible for the concentrations and the conditions (temperature and pressure)
arising during in practical applications, thus this direction can be omitted and
the kept direction can be considered as an irreversible step without inducing significant error.

Accordingly, in widely accepted combustion models irreversible steps are used, nevertheless, during model construction
we assume that the reaction steps are reversible and the omission one of the directions can be investigated a posteriori.
\begin{itemize}
\item
Hydrogen combustion:$^{\cite{healykalitanaulpetersenbourquecurran}}$
\begin{equation}\label{eq:R3}
\ce{H2O2 + O -> HO2 + OH}.
\end{equation}
\item
Carbon monoxide combustion:$^{\cite{healykalitanaulpetersenbourquecurran}}$
\begin{equation}\label{eq:R4}
\ce{HCO + HO2 -> CO2 + H + OH}, \ce{H2O2 + O -> HO2 + OH}.
\end{equation}
Further examples e.g. in methanol combustion ($^{\cite{arandachristensenalzuetaglarborggersencgaodmarshall}}$)
can also be found.
\torol{
or:$^{\cite{keromnesmetcalfeheuferdonohoedassungherzlernaumanngriebelmathieukrejcipetersenpitzcurran}}$
$$
\ce{2HCO -> 2 CO + H2} , \ce{HCO + HO2 -> CO2 + H + OH}.
 $$
\item
Methanol combustion:$^{\cite{arandachristensenalzuetaglarborggersencgaodmarshall}}$
\begin{align*}
&\ce{CH2OOH -> CH2O + OH},\\
&\ce{C2H4 + HOCH2CH2OO -> CH2O + CH2OH + CH3CHO},\\
&\ce{CH2O + HOCH2CH2OO -> CH2OH + CH2OOH + HCO},\\
&\ce{HO2 + HOCH2CH2OO -> CH2OH + CH2OOH + O2},\\
&\ce{CH2CHOH + O2 -> CH2O + HCO + OH}.
\end{align*}
Another pair of irreversible steps, namely
$$
\ce{NO2 -> NO2*}, \ce{2NO2* -> 2NO + O2}
$$
(valid for low concentrations of \ce{NO}) can be found in Ref. \cite{rasmussenwassarddamjohansenglarborg}.
}
\end{itemize}
\paragraph{\ce{X <=> X}, \ce{2X <=> 2X}, \ce{X + Y <=> X + Y}}
Certainly there is no reason to include reaction steps which do not affect concentrations
(i.e. make no macroscopic change) despite that they are taking place microscopically.
\paragraph{Mass conservation}
We discard steps like \ce{X <=> 2X}, \ce{X <=> X + Y}, as they obviously violate the law of mass conservation.
However, these kinds of steps may be quite useful in model construction in chemical kinetics.
E.g. the Lotka--Volterra reaction which was proposed
for \emph{approximately} describing oscillations in cold flames$^{\cite{frankkamenetskii}}$,
contains steps like \ce{X -> 2X} and \ce{Y -> 0}.
\paragraph{\ce{2X <=> 2Y}}
In models one can rarely see such a step.
However, studying the literature one can find similar steps: the self-reaction of peroxy radicals, like 
$$
\ce{2 CH3O2 <=> 2 CH3O} (+ \ce{O2}),
$$
or
$$
\ce{2 CH3CH2O2 <=> 2 CH3CH2O} (+ \ce{O2}),
$$
see Refs. \cite{horiecrowleymoortgat} and \cite{noellalconcelrobichaudokumurasander}.
One can omit the forming \ce{O2} in atmospheric (e.g. tropospheric) chemistry as it is a major
constituent of air and its concentration will be negligibly affected by this and similar transformations.

Both steps are of the form \ce{2X <=> 2Y} (+\ce{Z}) (and are mass conserving---if \ce{O2} is taken into consideration---as they fulfill the law of atomic balance).
\paragraph{Reactions passing all the criteria}
\paragraph{\fbox{$M=2$}}
The species, complexes and reaction steps are shown in the Table below.

\medskip
\begin{table}
\begin{center}
\begin{tabular}{|l|l|l|}
  \hline
  Species   & \ce{X}, \ce{Y}&$M = 2$ \\
  \hline
  Complexes &\ce{X}, \ce{Y}, \ce{2X}, \ce{2Y}, \ce{X + Y}&$\overline{N}(2) = 5$ \\
  \hline
  Reaction steps& \ce{X <=> Y}, \ce{X <=> 2 Y}, \ce{Y <=> 2 X}, &$\overline{R}(2) = 5$\\
   & \ce{2 X <=> X + Y}, \ce{2 Y <=> X + Y}& \\
  \hline
\end{tabular}
\caption {Species, complexes and reaction steps in the case $M=2$} \label{tab:mistwo}
\end{center}
\end{table}
\medskip
\paragraph{\fbox{$M=3$}}
The species, complexes and reaction steps are shown in the Table below.

\medskip
\begin{table}
\begin{center}
\begin{tabular}{|l|l|l|}
  \hline
  Species   & \ce{X}, \ce{Y}, \ce{Z}&$M = 3$ \\
  \hline
  Complexes &\ce{X}, \ce{Y}, \ce{Z}, \ce{2X}, \ce{2Y}, \ce{2Z}, &$\overline{N}(3) = 9$ \\
  &\ce{X + Y}, \ce{Y + Z}, \ce{Z + X}&\\
  \hline
  Reaction& \ce{X <=> Y}, \ce{X <=> Z}, \ce{Y <=> Z}, &$\overline{R}(3) = 24$\\
  steps&\ce{X <=> 2 Y}, \ce{X <=> 2 Z}, \ce{Y <=> 2 X}, & \\
  &\ce{Y <=> 2 Z}, \ce{Z <=> 2 X}, \ce{Z <=> 2 Y},&\\
  &\ce{X <=> Y + Z},\ce{Y <=> X + Z},&\\
  &\ce{Z <=> X + Y},\ce{2 X <=> X + Y},&\\
  & \ce{2 X <=> X + Z}, \ce{2 X <=> Y + Z},&\\
  &\ce{2 Y <=> X + Y}, \ce{2 Y <=> X + Z},&\\
  & \ce{2 Y <=> Y + Z},\ce{2 Z <=> X + Y},&\\
  & \ce{2 Z <=> X + Z}, \ce{2 Z <=> Y + Z},&\\
  &\ce{X + Y <=> X + Z}, \ce{X + Y <=> Y + Z},&\\
  & \ce{X + Z <=> Y + Z}&\\
  \hline
\end{tabular}
\caption {Species, complexes and reaction steps in the case $M=3$} \label{tab:misthree}
\end{center}
\end{table}

Now the interesting (from the combinatorial point of view) question arises:
What is the number of reactions fulfilling the requirements formulated above?
Beyond combinatorics, formulae to give those numbers is of practical relevance too:
it gives us a hint if it is possible to deal with all the systems constructed in this way within a tolerable time.
We enumerated the steps from \(M=2\) up to \(M=20\) an have found the following cardinalities:
$5,24,69,155,300,525,854\dots.$
How to learn if there is a certain regularity in the sequence?
The best way is to go to the Online Encyclopedia of Integer Sequences (\url{https://oeis.org/})  initiated by Neil James Alexander Sloane, and ask if it contains our sequence.
In this case the answer was yes, and the formula $(M-1)M(M^2+7M+2)/8$ is provided to give the number of reaction steps of the given type.
From the strict mathematical point of view now we had this statement as a conjecture, but it can rigorously be proved, as well.
\begin{Sta}
Suppose the number of species is \(M\in\N.\)
Then, the number of reversible, mass conserving reaction steps excluding steps of the form
$$
\ce{X <=> X}, \ce{2X <=> 2X} \mbox{ and } \ce{2X <=> 2Y}
$$ is
\begin{equation}\label{eq:sloane2}
\overline{R}(M):=(M-1)M(M^2+7M+2)/8, M=2,3,\dots.
\end{equation}
\end{Sta}
\Proof{
The table below shows how many steps of the different types we have.
\begin{table}
\begin{center}
\begin{tabular}{|c|c|r|r|r|}
\hline
Type of step&$M=$&2&3&4\\
\hline
\ce{X <=> Y}&$\frac{M(M-1)}{2}$&1&3&6\\
\hline
\ce{X <=> 2Y}&$M(M-1)$&2&6&12\\
\hline
\ce{2X <=> X + Y}&$M(M-1)$&2&6&12\\
\hline
\ce{X + Y <=> Z}&$\frac{M(M-1)}{2}(M-2)$&0&3&12\\
\hline
\ce{X + Y <=> 2Z}&$\frac{M(M-1)}{2}(M-2)$&0&3&12\\
\hline
\ce{X + Y <=> X + Z}&$\frac{M(M-1)}{2}(M-2)$&0&3&12\\
\hline
\ce{X + Y <=> Z + A}&$\frac{1}{2}\frac{M(M-1)}{2}\frac{(M-2)(M-3)}{2}$&0&0&3\\
\hline
Total&$\overline{R}(M)$&5&24&69\\
\hline
\end{tabular}
\caption{Number of reaction steps of different type}\label{tab:steps}
\end{center}
\end{table}
}
Below we will have other sequences of numbers to be studied.
\subsection{Reaction mechanisms}
A set of reaction steps is usually called a (kinetic reaction) \emph{mechanism}.
In formal reaction kinetics alternative names (complex chemical) \emph{reaction}, or \emph{reaction network} are also used.
\paragraph{The actual number of species is $M$}
We may start from three species and arrive at a mechanism with two species, as e.g.
$$
\ce{X <=> Y}, \ce{X <=> 2 Y}, \ce{Y <=> 2 X}
$$
may be obtained as one of the mechanisms with three species,
but it actually contains only two species.
It would be desirable to exclude such cases.

Suppose we have $M=3$ species, then 
the number of mechanisms \(\overline{R}(3)\) with not more than three species, and also
the number of mechanisms \(\tilde{R}(3)\) with exactly three species
 are shown in the Table below. 
\begin{table}\label{tab:less}
\begin{center}
\begin{tabular}{|l|r|r|r|r|r|r|r|}
\hline
Number of steps&1&2&3&4&5&6&$\ge6$\\
\hline
$M=3$   &9 & 246& 1994& 10611& 42501&134596&$\tilde{R}(3)$\\
$M\le 3$&24& 276& 2024& 10626& 42504&134596&$\binom{\overline{R}(3)}{R}$\\
\hline
\end{tabular}
\end{center}
\caption{Number of mechanisms containing different number of reaction steps}
\end{table}

In the case when we have at most three species it may happen that we only have 2,
thus the mechanisms with exactly $3$ species is 
$\tilde{R}(3)=\binom{\overline{R}(3)}{R}-\binom{\overline{R}(2)}{R}\binom{3}{2}$ 
because of the three species one can select any two in\(\binom{3}{2}\) different ways.
As $\overline{R}(2)=5,$ last (and following) elements in the two rows of the table coincide,
as one can only have 5 reaction steps with two species, see Table \ref{tab:mistwo}.
Let us formulate the corresponding---obvious---general statement.
\begin{Sta}
Suppose the number of species is \(M\in\N.\)
Then, the number of reversible, mass conserving mechanisms excluding steps of the form
\ce{X <=> X}, \ce{2X <=> 2X} and \ce{2X <=> 2Y} which do contain \(M\) species and consist of $R$
reactions is
\begin{align}
\tilde{R}(M)&:=
\binom{\overline{R}(M)}{R}
-\binom{\overline{R}(M-1)}{R}\binom{M}{M-1}.\label{eq:sloane32}
\end{align}
Consequently, in case $R>\overline{R}(M-1)$ one only has the first term.
\end{Sta}
\paragraph{Mass conservation}\label{mcreactions}
Even if the steps are mass conserving the mechanism may not be as the example 
\{\ce{X <=> Y},\ce{Y <=> 2X}\} shows.
To check this property is not a trivial problem, we do this using our program 
\mma{ReactionKinetics} described in Chapter 4 of Ref \cite{tothnagypapp}, 
where the reader can also find relevant references, as well.
\paragraph{Detailed balancing}
Following Section 7.8 of Ref. \cite{tothnagypapp} 
we review the history of detailed balancing shortly.

\begin{figure}
  \centering
  \includegraphics[width=300pt]{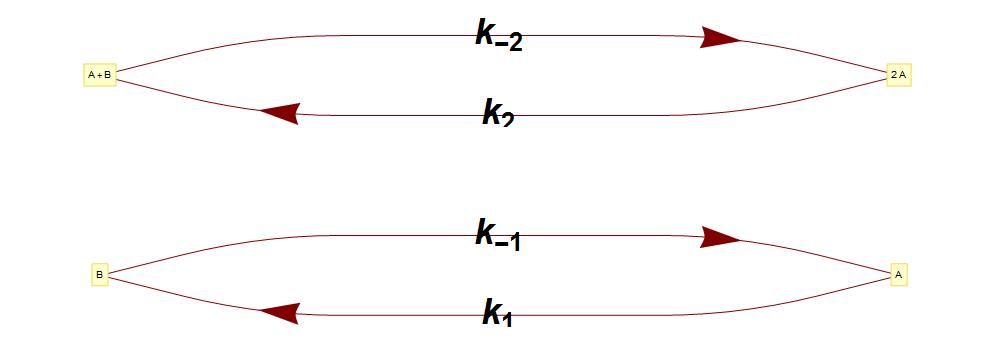}
  \caption{The Wegscheider mechanism}\label{fig:wegscheider}
\end{figure}

After such men as Maxwell and Boltzmann, and before Einstein
(see the references here: \url{https://en.wikipedia.org/wiki/Detailed_balance}),
at the beginning of the twentieth century, it was Wegscheider$^{\cite{wegscheider}}$ 
who constructed the reaction mechanism in Fig. \ref{fig:wegscheider}
to show that in some cases the existence of a positive stationary state 
alone does not imply the equality of all the individual forward and backward 
reaction rates in equilibrium: A relation (in this case $k_{-1} k_2 =k_{-2}k_1$) 
should hold between the reaction rate coefficients to ensure this. 
Equalities of this kind will be called (and later exactly defined) as
\emph{spanning forest conditions} below.
Let us emphasize that violation of this equality
between the reaction rate coefficients does not exclude the existence of a positive stationary state; 
it can be shown to exist and unique for all values of the reaction rate coefficients.
(Problem 7.12 of Ref. \cite{tothnagypapp} proves both statements.)

Here we mention that Figs. 1--3. show the Feinberg--Horn--Jackson graph of the mechanism
which is a directed graph with the complexes
as vertices and with the reaction step arrows as directed edges. 
The number of complexes is denoted by $N,$ the number of connected components of the 
Feinberg--Horn--Jackson graph is $L,$
whereas the number of independent reaction steps (the rank of $\gammab$) is $S.$

\begin{figure}
  \centering
  \includegraphics[width=300pt]{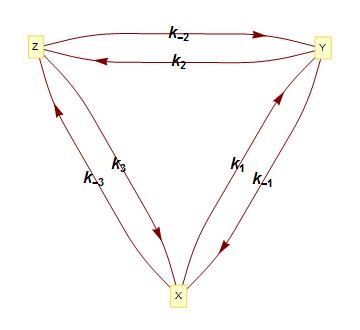}
  \caption{The reversible triangle reaction}\label{fig:revtri}
\end{figure}

A similar statement holds for the reversible triangle reaction in Fig. \ref{fig:revtri}.
The necessary and sufficient condition for the existence of such a positive stationary state 
for which all the reaction steps have the same rate in the forward and backward direction 
(a detailed balanced stationary state) is 
now that the product of the reaction rate coefficients is the same if taken in either direction:
$k_1k_2k_3=k_{-1}k_{-2}k_{-3}.$
Equalities of this kind will be called as \emph{circuit conditions} below.
Again, violation of this equality does not exclude the existence of a positive stationary state;
it can again be shown to exist and be unique for all values of the reaction rate coefficients.
(Problem 7.11 of the reference mentioned shows both statements.)

These examples are qualitatively different from, e.g., the simple one-step mechanism 
\ce{X + Y <=>[$k_1$][$k_{-1}$] Z} which has the same stationary reaction rate 
in both directions no matter what the values of the reaction rate coefficients are, 
see Problem 7.10.
To put it another way, this mechanism is \emph{unconditionally detailed balanced} while the
previous examples (being detailed balanced only if come equalities are fulfilled) were
\emph{conditionally detailed balanced}. To show a less trivial example, have a look at Fig.
\ref{fig:db}.

A quarter of a century after Wegscheider, the authors Fowler and Milne$^{\cite{fowlermilne}}$
formulated in a very vague form a general principle called the \emph{principle of detailed
balance} stating that in real thermodynamic equilibrium, all the subprocesses
should be in dynamic equilibrium
separately in such a way that they do not stop but proceed with the same velocity
in both directions.
Obviously, this also means that time is reversible at equilibrium;
that is why this property may also be called \emph{microscopic reversibility}, although it
may be appropriate to reserve this expression for a similar property of the stochastic
model (see Chap. 10 of Ref. \cite{tothnagypapp}).
A relatively complete summary of the early developments was given by Tolman$^{\cite{tolman}}$.
The modern formulation of the principle accepted by IUPAC$^{\cite{goldloeningmacnaughtshemi}}$
essentially means the same (given that the principle of charity is applied when reading):
“The principle of microscopic reversibility at equilibrium states that, in a
system at equilibrium, any molecular process and the reverse of that process
occur, on the average, at the same rate.”

Now we give a precise formulation of the concept in such a way that at a detailed
balanced stationary point (which can only exist in a reversible reaction), all pairs
of reaction–antireaction step pairs proceed with the same rate in both directions.

The reaction we study here consists of reversible pairs of reaction steps like
$$
\sum_{m=1}^M\alpha_{mr}\ce{X$_m$}\rightleftharpoons\sum_{m=1}^M\beta_{mr}\ce{X$_m$}
\quad(r=1,2,\dots,R)
$$
and their usual induced kinetic differential equations assuming mass action type kinetics
(and disregarding the change of temperature, pressure and reaction volume) is
$$
\dot{c}_m=\sum_{r=1}^P(\beta_{mr}-\alpha_{mr})
(k^+_r\prod_{p=1}^Mc_p^{\alpha_{pr}}-k^-_r\prod_{p=1}^Mc_p^{\beta_{pr}}),
$$
where $c_m(t):=[\ce{X}_m](t)$ is the concentration of species \ce{X$_m$}.
(Note that here $R$ denotes the half of the number of reaction steps, 
or the number of reversible pairs; 
we reserve the notation $P$ for the total number of reaction steps. 
Thus, in the case of reversible steps $P=2R.$)
Shortly,
$$
\dot{\cb}=\gammab (\kb^+\odot\cb^{\alphab}-\kb^-\odot\cb^{\betab})
$$
with
$\alphab:=(\alpha_{mr}), \betab:=(\beta_{mr}), \gammab:=\betab-\alphab,$
$\cb^{\alphab}:=\prod_{m=1}^{M}c_m^{\alpha_{mr}},$ 
and with the componentwise (or Schur) product $\odot$ of vectors.
Here the positive numbers $k^{\pm}_r$ are the \emph{reaction rate coefficients}, the vectors formed from them are $\kb^{\pm}.$
The reaction is detailed balanced at the positive stationary concentration $\cb_*$ if all the steps proceed with the same rate in both directions, or, to put it another way
$$
\text{DB: }
\gammab (\kb^+\odot\cb_*^{\alphab}-\kb^-\odot\cb_*^{\betab})=\nulb
\text{ implies }
$$
$$
\kb^+\odot\cb_*^{\alphab}=\kb^-\odot\cb_*^{\betab} \text{ or }
\gammab^\top\log(\cb_*)=\log(\Kb),\label{eq:forfredholm}
$$
where $\Kb:=\frac{\kb^+}{\kb^-}.$

Detailed balance may hold
\begin{itemize}
\item at any (positive) values of the reaction rate coefficients (unconditionally detailed balanced), or
\item only if the values of the rate coefficients fulfil certain conditions---e.g. 
circuit or spanning tree conditions---(conditionally detailed balanced).
\end{itemize}

What are the necessary and sufficient conditions of this property? 
First we give an algebraic characterization that can be proved using Fredholm's alternative theorem.
\begin{Thm}[See Refs. \cite{tothnagypapp,vladross}]\label{th:fredholm}
The reaction is detailed balanced, if and only if for all nonzero \emph{vector} solutions to the system of linear equations \(\gammab \ab=\nulb\) one has
\begin{equation}
\Kb^\ab=\oneb.\label{eq:db}
\end{equation}
\end{Thm} 
As the elements of \(\gammab\) are integers and the vectors \(\ab\) are solutions of 
a homogeneous linear equations, their coordinates can supposed to be integers.

Next we cite a pair of structural criteria showing what the reasons of detailed balancing are.
To formulate this, we need a few concept and also a few formal definitions.

\begin{Def}
The \emph{circuit conditions} are that the product of reaction rate coefficients along any set of independent cycles is the same in both directions.
\end{Def}
\begin{Def}
Let us take a spanning forest of the Feinberg--Horn--Jackson graph, and let the corresponding
reaction step vectors be $\gammab_{.,u}{\ }(u=1,2,\dots,N-L).$ 
Then, $\sum_{u}^{N-L}a_{u}\gammab_{.,u}=\nulb$ has $N-L-S$ independent solutions.
With these
$\prod \frac{k_{u}^+}{k_{u}^-}^{a_{u}}=1$  should hold: 
these are the \emph{spanning forest conditions}.
\end{Def}
Note that the number of the edges of the spanning \emph{tree} 
is \(L\) less than the number of its vertices, if again \(L\) is the number of
the connected components of the Feinberg--Horn--Jackson graph.
\begin{Thm}[Ref. \cite{feinbergdb}]\label{th:feinberg}
The mechanism is detailed balanced, if and only if the \emph{circuit conditions} and the \emph{spanning forest conditions} hold.
\end{Thm}

An application of Feinberg's theorem (and also the detailed description 
with examples of the concepts) can be found in Ref. \cite{nagytothion}.

\begin{Exa}
An unconditionally detailed balanced reaction can be seen in Figure \ref{fig:db}.
\begin{figure}
  \centering
  \includegraphics[width=300pt]{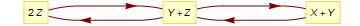}
  \caption{No matter what the values of the reaction rate coefficients are, 
  this is detailed balanced.}
  \label{fig:db}
\end{figure}
The reason is that both structural conditions are empty:
\begin{enumerate}
\item
It does not contain cycles.
\item
Its deficiency ($:=N-L-S=3-1-2$) is zero, thus no spanning forest conditions are to be considered.
\end{enumerate}
With the approach of applying Theorem \ref{th:fredholm} one sees that $\gammab\ab=\nulb$
has no nonzero solutions, or the kernel of the linear map $\gammab$ only contains the 
zero vector.
\end{Exa}
\begin{Exa}
A conditionally detailed balanced mechanism is shown in Figure \ref{fig:cdb}.
Here the spanning forest conditions (their number is $N-L-S=5-2-1=2$) are as follows.
\begin{equation}
k_{-1}^2k_3=k_{-3} k_1^2,{\ }
k_{-1} k_2=k_{-2}k_1.\label{eq:spanning}
\end{equation}
We also have the circuit conditions:
\begin{equation}
k_{-4} k_{-2} k_3=k_{-3} k_2 k_4.\label{eq:circuit}
\end{equation}
\begin{figure}
\centering
  \includegraphics[width=300pt]{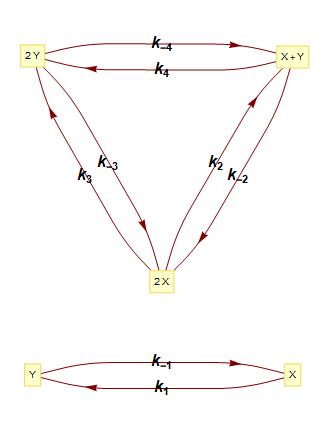}
  \caption{Detailed balanced, if the system of equalities \eqref{eq:spanning} and \eqref{eq:circuit} holds.}  \label{fig:cdb}
\end{figure}
\end{Exa}

Let us try to discuss unconditionally and conditionally detailed balanced reactions in a more systematic way.

If \fbox{$M=2, R=2$}, then the following three mass conserving mechanisms remain.
\begin{align}
&\ce{X <=>[$k_1$][$k_{-1}$] Y}, \ce{2 X <=>[$k_2$][$k_{-2}$] X + Y},&k_{-2}=\frac{k_{-1}k_2}{k_1}\label{eq:first}\\
&\ce{X <=>[$k_1$][$k_{-1}$] Y}, \ce{2 Y <=>[$k_2$][$k_{-2}$] X + Y},&k_{-2}=\frac{k_1k_2}{k_{-1}}\label{eq:second}\\
&\ce{2 X <=>[$k_1$][$k_{-1}$] X + Y}, \ce{2 Y <=>[$k_2$][$k_{-2}$] X + Y}&k_{-2}=\frac{k_1k_2}{k_{-1}}.\label{eq:third}
\end{align}
They are all conditionally detailed balanced.
(One can easily show that with $M=2$ and 2 reversible steps there are no unconditionally detailed balanced mechanisms.)
They contain no circles, thus only the \emph{spanning forest conditions} should hold, and these are shown in the second column above.

If \fbox{$M=3, R=2$}, there are
18 conditionally detailed balanced, 189 unconditionally detailed balanced mechanisms.
Let us calculate the number of conditionally and unconditionally detailed balanced reactions
for the cases $M=3,4.$ 
Then, the following table results, where \textcolor{brown}{MC} denotes the total number of mass conserving mechanisms, \textcolor{red}{UDB} is for unconditionally detailed balanced mechanisms, and CDB stands for conditionally detailed balanced mechanisms:
\textcolor{brown}{MC}=\textcolor{red}{UDB}+CDB

\begin{table}
\begin{center}\label{tab:udbcdb}
{\footnotesize
\begin{tabular}{|r|r|r|r|}
  \hline
   & $R=1$ & 2 & 3 \\
  \hline
  $M=3$ & 24 \textcolor{brown}{24}=\textcolor{red}{24}+0 & 276 \textcolor{brown}{207}=\textcolor{red}{189}+18 & 2024 \textcolor{brown}{602}=\textcolor{red}{0}+602 \\
  4 & 69 \textcolor{brown}{69}=\textcolor{red}{69}+0 & 2346 \textcolor{brown}{2064}=\textcolor{red}{2004}+60 & 52394 \textcolor{brown}{30768}=\textcolor{red}{2478}+5898 \\
  \hline
\end{tabular}
}
\caption{The number of not mass conserving and of mass conserving and 
either unconditionally detailed balanced or conditionally detailed balanced mechanisms}
\end{center}
\end{table}

How to use detailed balance in our calculations? 
We have two choices.
One can either fit using independent reaction rate coefficients and check detailed balance later, 
or use detailed balance as a condition when fitting. 
In our calculations we have used the second alternative.
\torol{
\paragraph{Fitting vector valued functions}\label{subsub:vector}
Suppose we are given a function (procedure, subroutine, application) like \mma{NonlinearModelFit}---of the Wolfram language---that is able to find the parameters of a scalar valued function, nonlinear in the parameters.
How to use it for vector valued functions?

Introduce
$$
z(1,t_i):=x(t_i)\quad z(2,t_i):=y(t_i)
$$
in case $x$ and $y$ are measured at times $t_i.$
Thus, instead of two scalar valued functions (with a scalar independent variable)
we arrived at a single scalar valued vector variable function.
We have used this trick in the calculations.

Note that the first independent variable of $z$ can only take the values 1 or 2, and ("therefore")
continuity in this variable is not assumed or fulfilled.
}
Let us mention that our calculations heavily rely on the program package \mma{ReactionKinetics} 
written in Wolfram language (Mathematica) and downloadable from \url{extras.springer.com} 
using the ISBN number 978-1-4939-8641-5.
This package is aimed at helping the chemist to do many kinds of symbolic and 
numerical investigations of reaction mechanisms
including solving the induced kinetic differential equations or simulating the usual stochastic model, but excluding parameter estimation.
The codes written to the present paper will be provided to the reader upon request.
\section{APPLICATIONS OF THE METHOD}
Although to illustrate the method we use toy models, we think the results are promising.
Two applications of the method will be shown.
\begin{enumerate}
\item
Of the three reversible mechanisms constructed with $M=2$ species and $R=2$ 
pairs of reaction steps above we simulated data 
for one of them, and the program fitted all three candidate mechanisms
and identified the best of them.
As they are inherently conditionally detailed balanced,
detailed balance was used as a constraint, a prescribed relationship 
between the reaction rate coefficients during fitting.
\item
\emph{Real experimental data:} reliable old measurements$^{\cite{raczgyarmatitoth}}$ 
on the transport of salicylic acid with or without different additives under different 
pH values could be explained by the simple consecutive mechanism: \ce{X -> Y -> Z}.
Note that there are no chemical reactions involved, and
the same species in different compartments are denoted as different species and their transport is described as formal reactions.

Now we tried to fit all the six models \ce{X -> Y -> Z}, \ce{X -> Z -> Y}, etc.
generated by permuting the order of the compartments.
\end{enumerate}
\subsection{Simulated data in the $M=2,R=2$ case}
\begin{figure}
  \centering
  \includegraphics[width=150pt]{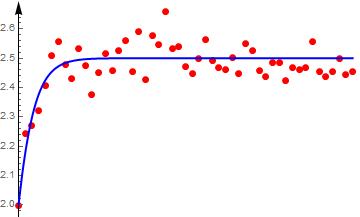}$\quad$\includegraphics[width=150pt]{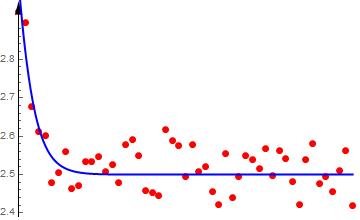}
  \caption{The original model: \ce{X <=>[$k_1$][$k_{-1}$] Y}, \ce{2 X <=>[$k_2$][$k_{-2}$] X + Y}}\label{fig:orig}
\end{figure}

We solved the deterministic model
\begin{align*}
&x'(t)=\frac{k_{-1} k_1}{k_2}x(t) y(t)-k_2 x(t)^2-k_1 x(t)+k_{-1} y(t),\\
&y'(t)=-\frac{k_{-1} k_1}{k_2}x(t) y(t)+k_2x(t)^2+k_1 x(t)-k_{-1} y(t)
\end{align*}
of the reaction
\ce{X <=>[$0.1$][$0.1$] Y}, \ce{2 X <=>[$1$][$k_{-2}$] X + Y}
with the given reaction rate coefficients and $k_{-2}$ calculated from the condition of detailed balance
\eqref{eq:first}.
The initial concentrations were $x(0)=2, y(0)=3.$
Then, we took a sample from the concentrations at equidistant discrete times with a sampling step size 0.1 and added a 2\% relative error.
Next, we tried to fit the model of the reaction from which we started.
We have chosen the initial estimates of the parameters to be $0.5, 0.5, 0.5.$
the estimated parameters became
$$
k_1 = 0.0971753, k_{-1} = 0.0971771, k_2 = 0.999665.
$$
Figure \ref{fig:orig} suggests a good fitting.

\begin{figure}
  \centering
  \includegraphics[width=150pt]{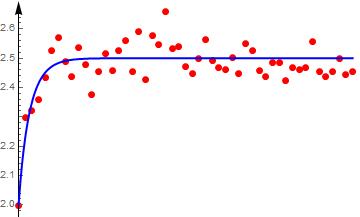}$\quad$\includegraphics[width=150pt]{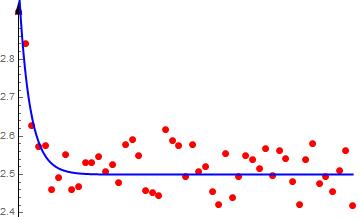}
  \caption{A slightly worse model: \ce{2 X <=> X + Y}, \ce{2 Y <=> X + Y}}\label{fig:third}
\end{figure}

In the second case everything remained unchanged except that the reaction to be fitted was
\ce{2 X <=> X + Y}, \ce{2 Y <=> X + Y}.
Note that condition of detailed balance is different from the first case again, see \eqref{eq:third}.
With the same data neither the numerical result, nor the fitting seems to be too bad (although slightly worse:
$$
k_1 = 0.0741455, k_{-1} = 0.0741468, k_2 = 1.01966.
$$
see also Figure \ref{fig:third}.

\begin{figure}
  \centering
  \includegraphics[width=150pt]{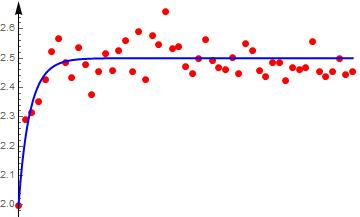}$\quad$\includegraphics[width=150pt]{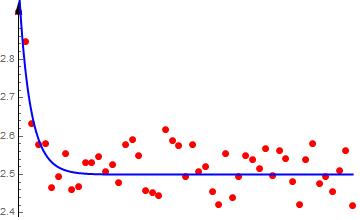}
  \caption{An unacceptable model: \ce{X <=> Y}, \ce{2 Y <=> X + Y}}\label{fig:second}
\end{figure}

Finally, the reaction to be fitted was
\ce{X <=> Y}, \ce{2 Y <=> X + Y}.
Note that condition of detailed balance is different in this case, see \eqref{eq:second}.
Again, with the same data the fitting does not seem to be too bad, but the estimated parameter values are unacceptable:
$$
k_1 = -0.357567, k_{-1} = -0.357585, k_2 = 1.16416.
$$
see also Figure \ref{fig:second}.
(One can force the program to provide positive reaction rate coefficients only,
but in the first two cases this was not needed.)

\subsection{Salicylic acid transport with $M=3,R=2$ case}
The concentrations of salicylic acid were measured in two compartments \ce{X} (gastric fluid) and
\ce{Z} (intestine fluid),
assuming that the process is described with the simple consecutive reaction \ce{X -> Y -> Z}.
The concentrations measured in the different compartments can be seen in the Table.

\medskip
\begin{tabular}{|r|r|r|r|}
  \hline
  $t_i$ (hours)& $x(t_i)$ (M/l) & $t_i$ (hours) & $z(t_i)$ (M/l) \\
  \hline
  1  & 0.01579 & 0 & 0       \\
  2  & 0.01429 & 1 & 0.0003  \\
  3  & 0.01327 & 2 & 0.000614\\
  4  & 0.01230 & 3 & 0.000917\\
  5  & 0.01148 & 4 & 0.00143 \\
  6  & 0.01066 & 5 & 0.00201 \\
  7  & 0.00988 & 6 & 0.00269 \\
  8  & 0.00912 & 7 & 0.00338 \\
  9  & 0.00851 & 8 & 0.00402 \\
  10 & 0.00791 & 9 & 0.00473 \\
  \hline
\end{tabular}
\medskip

\begin{figure}
  \centering
  \includegraphics[width=140pt]{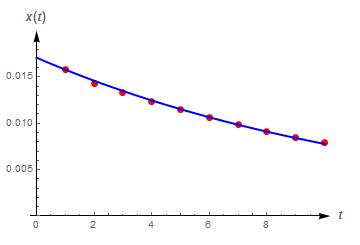}{\ }\includegraphics[width=140pt]{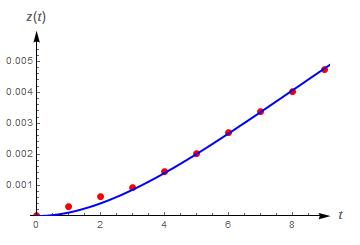}
  \caption{Fitting a bad model to real life data}\label{fig:salicyl}
\end{figure}

Fitting the mechanism \ce{X ->[$k_1$] Y ->[$k_2$] Z} gives nice results: starting from the initial estimates
$(1,1)$ for the reaction rate coefficients (more precisely, transport coefficients)
gives the results with small standard error
$\{k_1 -> 0.0786482\pm0.000812339, k_2 -> 0.181337\pm0.00527615\},$
a "good" correlation matrix
\[
\left[\begin{array}{rr}1.& -0.499719\\-0.499719&1.\end{array}\right].
\]
The fitting can also be seen to be good, see Fig. \ref{fig:salicyl}.

\begin{figure}
  \centering
  \includegraphics[width=140pt]{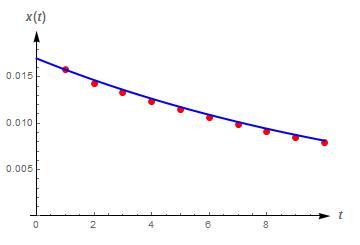}{\ }\includegraphics[width=140pt]{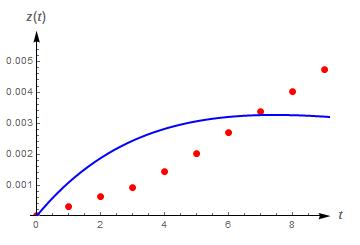}
  \caption{Fitting to real life data}\label{fig:salicylbad}
\end{figure}

However, if one wants to fit erroneously the model \ce{X -> Z -> Y}, then the \ce{Z} species fits badly,
as seen in Fig. \ref{fig:salicylbad}. Let us remark,
that the standard errors are also larger than previously,
which is not too interesting if one keeps in mind that the
transport coefficients of the bad model are meaningless.

\section{CONCLUSIONS, OUTLOOK}
At the end we make a few remarks on the restrictions used above, and also, on the final use of the results obtained from the type of calculations we have proposed.
\subsection{Restrictions}
We have shown two sets of reactions, 
but fitting to the given problem one can formulate specific sets of restrictions taking into 
consideration the chemical nature of the problem and trying to reduce the number of 
possible candidate models.
\paragraph{Mass action type kinetics}
We have almost always used it: especially in connection with detailed balance.
\paragraph{Complex balancing}
It is quite natural to ask why do we not select \emph{complex balanced reactions}, because this concept turned out to be more fundamental during the development of formal reaction kinetics$^{\cite{feinbergbook}}$.
Because in the case of mass action type kinetics we have a nice (structural) necessary and sufficient condition (the reaction should be weakly reversible and should have a zero deficiency $\delta:=N-L-S=0,$) this property can also be simply tested.
(As far as the concept goes, but one may have combinatorial blow up again.)

Detailed balancing and complex balancing are not so far from each other as they seem to be$^{\cite{dickensteinmillan,gorbanyablonsky,schaftraojayawardhana,szederkenyihangos}}$.
\paragraph{Atomic structure}
In this paper we only used chemical species without atomic structures.
From the knowledge of the atomic structures further restrictions are derived: it is natural to assume the law of atomic balance.

Combinatorial blow up.
Still, it will work for systems with not so many species.
Many things should be calculated once, it is only the fitting what takes time.

Our immediate goal is to extend and apply the method to multiple real life data sets.
\subsection{Final use of the the results}
Present the best fitting reaction using AIC (Akaike Information Criterion),
BIC (Bayesian Information Criterion) etc. to the chemist for further investigations and interpretation.
Here we do not want to spend to much time with the detailed statistical analysis of simulated data,
see the following books on this topics: \cite{bard,seberwild,turanyitomlin,weise}.

The method presented here can also applied as a special kind of lumping,$^{\cite{tothlirabitztomlin}}$.
Suppose a big kinetic model is given and we construct such small models in the way described above
which are able to reproduce the measured or simulated data obtained from the big one.

Carefully designed set of restrictions may decrease the number of candidate reactions,
and increasing capacity of computers can do away with larger numbers, too.

Our calculations and data can be requested from the authors.
\section*{NOTATIONS}
\begin{table}[!h]
  \centering
  \begin{tabular}{|l|l|l|l|}
     \hline
     Classical & Notation in & Meaning & Unit? \\
     notation&formal kinetics&Number of&\\
     \hline
      & L & linkage classes &  \\
      & M & species &  \\
      & N & complexes &  \\
      & P & reaction steps&\\ 
      & R & reversible pairs of reactions &  \\
      & S & independent reactions &  \\
     \hline
   \end{tabular}
  \caption{Notations}\label{tab:notations}
\end{table}

\section*{ACKNOWLEDGEMENTS}
Irene Otero-Muras raised complex balance as possible property to use for restriction.
Judit Z\'ador has supplied us with (real life chemical) references continuously.
The present work has been supported by the National Research, Development and Innovation Office
Hungary (SNN 125739).
\section*{ORCID}
\textit{Tibor Nagy} \url{https://orcid.org/0000-0002-1412-3007}

\noindent\textit{J\'anos T\'oth} \url{https://orcid.org/0000-0003-3065-5596}

\noindent\textit{Tam\'as Ladics} \url{https://orcid.org/}

\bibliography{NagyTothLadics}
\bibliographystyle{plain}
\end{document}